\documentstyle[12pt,amssymb]{amsart}

\begin{document}

\newtheorem{Theorem}{Theorem}
\newtheorem{Lemma}{Lemma}
\newtheorem{Conjecture}{Conjecture}
\newtheorem{Corollary}{Corollary}
\newtheorem{Definition}{Definition}
\newtheorem{Proposition}{Proposition}
\let\emptyset\varnothing

\newcommand{\D}{\bf D}
\newcommand{\T}{\bf T}

\def\supp{\mathrm{supp}\,}

\title[Douglas Algebras With No Maximal Sub- or Minimal Superalgebra]
{Douglas Algebras With No Maximal Subalgebra and No Minimal
Superalgebra}
\author{Carroll Guillory}
\address{University of Southwestern Louisiana\\
Lafayette, Louisiana 70504}
\curraddr{MSRI, 1000 Centennial Drive, Berkeley, CA 94720}
\email{cjg@@usl.edu}
\thanks{Research at MSRI is partially supported by NSF grant DMS-9022140.}

\maketitle

\begin{abstract}
We give several examples of Douglas Algebras that do not have any maximal
subalgebra. We find a condition on these algebras that guarantees that
some do not have any minimal superalgebra.
We also show that if $A$ is the only maximal subalgebra of  a Douglas
algebra $B$, then the algebra $A$ does not have any maximal subalgebra.
\end{abstract}

\section{Introduction}
Let $\D$ denote the open unit disk in the complex plane and $\T$ the unit
circle.
By $L^{\infty}$ we mean the space of essentially bounded measurable
functions on $\T$ with respect
to the normalized Lebesgue measure. We denote by $H^{\infty}$ the space of
all bounded analytic functions in $\D$.
Via identification with boundary functions, $H^{\infty}$ can be considered
as a uniformly
closed subalgebra of $L^{\infty}$. Any uniformly closed subalgebra $B$ strictly
between $L^{\infty}$ and $H^{\infty}$ is called a Douglas algebra.
 $M(B)$ will denote the maximal ideal space of a Douglas algebra $B$.
 If we let $X=M(L^{\infty})$, we can identify $L^{\infty}$
with $C(X)$, the algebra of continuous functions on $X$.
 If $C$ is the set of all continuous functions on $\T$, we set
$$H^{\infty}+C=\{h+g : g \in C , \, h \in H^{\infty} \}$$
$H^{\infty}+C$ then becomes the smallest Douglas algebra
containing $H^{\infty}$ properly.
The function
$$q(z)=\prod_{n=1}^{\infty}\frac{\left|z_n\right|}{z_n}\frac{z-z_n}{1-\bar
z_nz}$$
is called a Blaschke product if $\sum_{n=1}^{\infty}\left(1-|z_n|\right)$
converges. The set $\{z_n\}$ is called the zero set of $q$ in $\D$. Here
$\frac{|z_n|}{z_n}=1$ is understood whenever $z_n=0$. We call $q$ an
interpolating
Blaschke product if
$$\inf_{n}\prod_{m\neq n}\left|\frac{z_m-z_n}{1-\bar z_nz_m}\right|
>0 \cdot$$
An interpolating Blaschke product $q$ is called sparse (or thin) if
$$\lim_{n\to \infty} \prod_{m\neq
n}\left|\frac{z_m-z_n}{1-\bar z_nz_m}\right|=1 \cdot$$
The set
$$Z(q)=\left\{ x \in M(H^{\infty})\setminus \D : q(x)=0 \right\} $$
is called the zero set of $b$ in $M(H^{\infty}+C)$.
Any function $h$ in $H^{\infty}$ with $|h|=1$ almost everywhere on $\T$ is
called an inner
function. Since $|q|=1$ for any Blaschke product, Blaschke products are
inner functions.
Let
$$QC= \left(H^{\infty}+C\right)\cap \overline{\left(H^{\infty}+C\right)}$$
and, for $x \in M(H^{\infty}+C)$, set
$$Q_x=\left\{y \in M(L^{\infty}) : f(x)=f(y) \,\mathrm{for \,\,all}\,f\in
QC\right\}\cdot$$
 $Q_x$ is called the QC-level set for $x$. For $x\in M(H^{\infty}+C)$, we
denote
by $u_x$ the representing measure for $x$ and its support set by
$\supp u_x$. By
$H^{\infty}[\bar{q}]$ we mean the Douglas algebra generated by
$H^{\infty}$ and the complex
conjugate of the function $q$. Since $X$ is the Shilov boundary for every
Douglas algebra,
a closed set $E$ contained in $X$ is called a peak set for a Douglas
algebra $B$ if there is
a function in $B$ with $f=1$ on $E$ and $|f| < 1$ on $X\setminus E$. A
closed set $E$
is a weak peak set for $B$ if $E$ is the intersection of a family of peak sets.
If the set $E$ is a weak peak set for $H^{\infty}$ and we define
$$H^{\infty}_E=\left\{f\in L^{\infty} : \left. f\right|_E \in
\left.H^{\infty}\right|_E\right\},$$
then $H^{\infty}_E$ is a Douglas algebra. For a Douglas algebra $B$, $B_E$
is similarly defined.
 A closed set $E$ contained in $X$ is called the essential set for $B$, denoted
by $\mathrm{ess}(B)$, if $E$ is the smallest set in $X$ with the property
that for $f$ in
$L^{\infty}$ with $f=0$ on $E$, then $f$ is in $B$.

For an interpolating Blaschke product $q$ we put $N(\bar{q})$ the closure of
$$\bigcup \left\{ \supp u_x : x \in
M(H^{\infty}+C)\,\mathrm{and}\,|q(x)|<1\right\}\cdot$$
$N(\bar{q})$ is a weak peak set for $H^{\infty}$ and is referred to as
the nonanalytic
points of $q$. By $N_0(\bar{q})$ we denote the closure of
$$\bigcup \left\{ \supp u_x : x \in Z(q) \right\} .$$
For an $x \in M(H^{\infty})$ we let
$$E_x=\left\{ y \in M(H^{\infty}) : \supp u_y=\supp u_x\right\}$$
and call $E_x$ the level set of $x$.
Since the sets $\supp u_x$ and
$N(\bar{q})$ are weak peak sets for $H^{\infty}$, both
$H^{\infty}_{\supp u_x}$ and $H^{\infty}_{N(\bar{b})}$ are
Douglas algebras.
 For any interpolating Blaschke product $q$ we set
$$A=\bigcap_{x\in M(H^{\infty}+C)\setminus
M(H^{\infty}[\bar{q}])}H^{\infty}_{\supp u_x}$$
and
$$A_0=\bigcap_{y\in Z(q)}H^{\infty}_{\supp u_y}.$$

It is easy to see that $A\subseteq A_0$ and it was
shown in [11] that $A=H^{\infty}_{N(q)}$.
For $x$ and $y$ in $M(H^{\infty})$,
 the pseudo-hyperbolic distance is defined by
$$\rho (x,y)=\sup \left\{|h(x)| : |h| \leq 1, h \in H^{\infty}, h(y)=0
\right\} \cdot$$
For $x$ and $y$ in $\D$ we have$$\rho
(x,y)=\left|\frac{x-y}{1-\bar{y}x}\right| \cdot$$
For any $x \in M(H^{\infty})$, we define the Gleason part of $x$ by
$$P_x=\{y\in M(H^{\infty}) : \rho (x,y) <1 \}\cdot$$

If $P_x\neq \{x\}$, then $x$ is said to be a nontrivial point. We denote by $G$
the set of nontrivial points of $M(H^{\infty}+C)$, and for a Douglas
algebra $B$, we set
$$G_B=G\cap\left(M(H^{\infty}+C)\setminus M(B)\right) \cdot$$
A point $x$ in $G_B$ is called a minimal support point of $G_B$ (or simply
a minimal
support point of B) if there are no $y \in G_B$ such that
$\supp u_y \subseteq \supp u_x$. The set $\supp u_x$
is called a minimal support set for $B$. For Douglas algebras $B$ and $B_0$
with
$B_0 \subseteq B$ we let $\Omega (B,B_0)$ be all interpolating Blaschke
products $q$
such that $\bar{q} \in B$ but $\bar{q} \not\in B_0$.

We denote by $\Omega (B)$
the set of all interpolating Blaschke
 products $q$
with $\bar{q} \in B$.
Let $B$ be a Douglas algebra.
The Bourgain algebra $B_b$ of $B$ relative
to $L^{\infty}$ is the set of those elements
of $L^{\infty}$, $f$, such that $\| ff_n+B\|_{\infty} \to 0$
for every sequence $\{f_n\}$ in $B$ with $f_n \to 0$ weakly.
The minimal envelop $B_m$ of a Douglas algebra $B$
is defined to be the smallest Douglas algebra which contains all minimal
superalgebra of $B$. An algebra $A$ is called a minimal superalgebra of $B$
if for all
 $x,y \in M(B)\setminus M(A)$, $x\neq y$ implies
$\supp u_x=\supp u_y$.

A major part of this paper was done while the author was at the
Mathematical Sciences Research Institute. The author thanks the
Institute for its support and also thanks Pamela Gorkin for the many valuable
discussions.

\section{A Condition for Douglas Algebras to Have Equal Essential Sets}
Consider the Douglas algebras $A$ and $A_0$ defined above. In [11]
some conditions were given when $A \subsetneq A_0$ but yet
$\mathrm{ess}(A)=\mathrm{ess}(A_0)$. This happened because
$\mathrm{ess}(A)=N(\bar{q})$ and $\mathrm{ess}(A_0)=N_0(\bar{q})$
(this is not hard to show). Theorem 1 of [11] gives conditions when
$\mathrm{ess}(A)\neq \mathrm{ess}(A_0)$. The conditions found in Theorem 5
of [11] are far more complicated then those found in Theorem 1
below.
Yet $\mathrm{ess}(A)=\mathrm{ess}(A_0)$ in that Theorem (Theorem 5 of
[11]) and also satisfies the condition in Theorem 1.

\begin{Theorem}
Let $B_0$ be a subalgebra of a Douglas algebra $B$ with
$$\mathrm{ess}(B_0)\neq X.$$
If for every $x\in M(B_0)\setminus M(B)$ we have
$
\mathrm{ess}(H^{\infty}_{\mathrm{\supp
}u_x})=\mathrm{ess}(B_{\supp u_x}),
$
then $\mathrm{ess}(B)=\mathrm{ess}(B_0)$.
\end{Theorem}

{\it Proof:} 
We note that $\mathrm{ess}(H^{\infty}_{\supp u_x})=\supp u_x$.
Hence if $\mathrm{ess}(H^{\infty}_{\supp u_x})$ is contained in
$\mathrm{ess}(B)$
 for every $x\in M(B_0)\setminus M(B)$, then
$\supp u_y \subset \mathrm{ess}(B)$ for every
$y\in M(B_0)$ and so we get $\mathrm{ess}(B_0)\subset \mathrm{ess}(B)$.
Since $B_0\subset B$ we have that $\mathrm{ess}(B)\subset
\mathrm{ess}(B_0)$, and we get
$\mathrm{ess}(B)=\mathrm{ess}(B_0)$.

\begin{Corollary} Let $B_0$ be a maximal subalgebra of a Douglas algebra $B$.
Then $\mathrm{ess}(B_0)=\mathrm{ess}(B)$.
\end{Corollary}

Since $M(B_0)=M(B)\cup E_x$ for some
$x\in M(B_0)\setminus M(B)$, we have that if $z$ and $y$ are in
$M(B_0)\setminus M(B)$,
then $\supp u_y=\supp u_x=\supp u_z$. Now
$ess(B_{\supp u_x})=ess(H^{\infty}_{\supp u_x})$ since the set
$$\bigcup \left\{ \supp u_y : y\in M(B)\cap
M(H^{\infty}_{\supp u_x})\right\}$$
is dense in $\supp u_x$ (because $x$ is a minimal support point of $B$).
Thus $\mathrm{ess}(B)=\mathrm{ess}(B_0)$.

\begin{Corollary}
Let $A$ be any Douglas algebra and $q$ be an
interpolating Blaschke
product with $\bar{q} \not\in A$. Then
$\mathrm{ess}(A)=\mathrm{ess}(A[\bar{q}])$.
\end{Corollary}

{\it Proof:} Let $B=A[\bar{q}]$ and let $x\in M(A)\setminus M(B)$.
We need the following two facts:

\begin{enumerate}

\item[i)] $M(A)=\left\{x\in M(H^{\infty}+C) : A_{\supp u_x} =
 H^{\infty}_{\supp u_x}\right\}$

and

\item[ii)]  If $A\subset B$, then $A_{\supp u_x}\subset
B_{\supp u_x}$.

\end{enumerate}

We use   i) and ii) to show that
$B_{\supp u_x}=A_{\supp u_x}[\bar{b}]$.
Since $A\subset B$ by ii)
$M(B_{\supp u_x})\subset M(A_{\supp u_x})$.
By the Chang--Marshall Theorem [1,16] it suffices to show that
$M(B_{\supp u_x})=M(A_{\supp u_x}[\bar{q}])$.
Let $y\in M(B_{\supp u_x})$.
Since $M(B_{\supp u_x}) \subset M(A_{\supp u_x})$
 we have that $y\in M(A_{\supp u_x})$ and $|q(y)|=1$
(since $M(B)=\{ y\in M(A) : |q(y)|=1\}$). Hence $y\in
M(A_{\supp u_x}[\bar{q}])$.
 Therefore $M(B_{\supp u_x})\subset
M(A_{\supp u_x}[\bar{q}])$.

Now suppose $y \not\in M(B_{\supp u_x})$.
If $y\not\in M(A_{\supp u_x})$,
then $y\not\in A_{\supp u_x}[\bar{b}]$ and we have nothing to prove.
So we can assume that $y\in M(A_{\supp u_x})$. Now since
$y\not\in M(B_{\supp u_x})$ implies that $|q(y)|<1$ and
$y\in M(A_{\supp u_x})$.  Hence $y\not\in
M(A_{\supp u_x}[\bar{q}])$.
 Thus $M(A_{\supp u_x}[\bar{q}]) \subset M(B_{\supp u_x})$.
By the Chang--Marshall Theorem we have that
$M(A_{\supp u_x}[\bar{q}])=M(B_{\supp u_x})$.
Hence $A_{\supp u_x}[\bar{q}]=B_{\supp u_x}$.
Thus we have
\begin{eqnarray*}
\mathrm{ess}(B_{\supp u_x}) & = & \mathrm{ess}\left(
A_{\supp u_x}[\bar{q}]\right) \\
  & = & \mathrm{ess}\left(H^{\infty}_{\supp u_x}[\bar{q}]\right)
\quad \mathrm{by\ i)}
\end{eqnarray*}
So it suffices to show that
$$\mathrm{ess}\left(H^{\infty}_{\supp u_x}[\bar{q}]\right)=
\mathrm{ess}\left( H^{\infty}_{\supp u_x}\right).$$
To do this set $B_1=H^{\infty}_{\supp u_x}$ and
$$E=\bigcup \left\{ \supp u_y :
\supp u_y \subset\supp u_x,\,|q(y)|=1\right\}.$$
Assume that $\bar{E}$, the closure of $E$,
is properly contained in $\supp u_x$. Put
$$B_2=\bar H^{\infty}_E=
\overline{\left\{f\in L^{\infty} : \left. f\right|_E\in \left.
H^{\infty}\right|_E\right\}}.$$
By [2, page 39]
$M(B_2)=\left\{ m\in M(H^{\infty}) : \supp u_m \subseteq
\bar{E}\right\}
\cup M(L^{\infty})$.
Since $\bar{E} \subsetneq \supp u_x$ we have that
$B_1 \subsetneq B_2$. Therefore $M(B_2)
\subsetneq M(B_1)$ and so there is a nontrivial point
 $x_0\in M(B_1)\setminus M(B_2)$ [5, Proposition 4.1]
such that (a) $\supp u_{x_0} \subsetneq  \supp u_x$,
(b) $\supp u_{x_0} \not\subset \bar{E}$
[otherwise $x_0\in M(B_2)$], (c) $|q(x_0)|<1$ and (d)
$\supp u_{x_0}\cap E =\emptyset$. In (a)
$\supp u_{x_0}\subsetneq \supp u_x$ for if
$\supp u_{x_0} =\supp u_x$ for all $x_0\in M(B_1)\setminus
M(B_2)$,
then $B_1$ is a maximal subalgebra of $B_2$, so by Corollary 1 we have
$\bar{E}=\supp u_x=\mathrm{ess}(B_1)$ (in fact $B_2$
has no maximal subalgebras). The fact that $|q(x_0)|<1$ implies that
$\bar{q}
\not\in H^{\infty}_{\supp u_{x_0}}$.

Therefore there is a $y_0\in Z(q)$ such that
$\supp u_{y_0}\subseteq\supp u_{x_0}$.
 By Theorem 2 of [10]
 there is a $z_0\in Z(q)$
such that $\supp u_{z_0}$ is a minimal support
 set for $H^{\infty}[\bar{q}]$
that is contained in $\supp u_{x_0}\subseteq \supp u_x$.
Since $q$ is an interpolating Blaschke product $\supp u_{z_0}$ is
not trivial.
By [5, Theorem 4.2] there is $m\in M(H^{\infty}+C)$ so that
$\supp u_m$
 is nontrivial and
$\supp u_m\subsetneq \supp u_{z_0}\subseteq
\supp u_{x_0}$.
Since $\supp u_{z_0}$ is a minimal support set for
$H^{\infty}[\bar{q}]$
we have that $|q(m)|=1$. Thus $\supp u_{x_0}\cap E\neq \emptyset$.
This contradicts (d). So $\bar{E}=\supp u_x$ and
$\mathrm{ess}(H^{\infty}_{\supp u_x})=
\mathrm{ess}(H^{\infty}_{\supp u_x[\bar{q}]})$ and we get
$\mathrm{ess}(B_{\supp u_x})=\mathrm{ess}(H^{\infty}_{\supp u_
x})$
for every $x\in M(A)\setminus M(B)$, so Corollary 2 follows from our Theorem.

 We mention here that Corollary 2 was proved in [20, Theorem 2]
by another method.

There are algebras $B_0$ and $B$ that satisfy the
hypothesis of Theorem 1 and are not of the form $B=B_0[\bar{q}]$ for
any interpolating
Blaschke product (if $B_0\subseteq B$). To see this let $\Gamma$ be the
collection
of sparse Blaschke products and $B$ be the Douglas algebra
$[H^{\infty} : \bar{q};\,q\in \Gamma]$. Let $q_0$ be any element in
$\Gamma$ and put
$B_0=H^{\infty}[\bar{q}]$. Then $B_0\subset B$. By a Theorem of
Hedenmaln's [13]
we have that if $b$ is a Blaschke product such that $\bar{b}\in B$,
then $b=b_1\cdots b_n$ where each $b_i,\,\,i=1,\dots ,n$, is a
sparse Blaschke product. Hence if $x\in M(B_0)\setminus M(B)$, then $x$ is
the zero of some
sparse Blaschke product. So $x$ is a minimal support point of
$B$ for every $x\in M(B_0)\setminus M(B)$. This implies that
$\mathrm{ess}(B_{\supp u_x})=\mathrm{ess}(H^{\infty}_{\supp u_
x})$.
So, by Theorem 1, we have $\mathrm{ess}(B)=\mathrm{ess}(B_0)$
(Theorem 2 below shows that $H^{\infty}_{\supp u_x}$ is a maximal
subalgebra of
$B_{\supp u_x}$.) Now suppose there is a Blaschke product $q\in
\Omega (B,B_0)$
with $B=B_0[\bar{q}]$. Again by Hedenmaln's Theorem we have
$q=q_1\cdots q_n$ with each
$q_i$ a sparse Blaschke product. Let $Q$ be any infinite sparse Blaschke
product such that
$|Q|=1$ on $\bigcup _{x\in Z(q)}P_x$. Then there is a $m\in M(H^{\infty}+C)$
such that $Q(m)=0$ but $m\not\in \bigcup_{x\in Z(q)}P_x$.
Thus $|q(m)|=1$ and so we get that $m\in M(B_0[\bar{q}])$.
Thus $\bar{Q}\not\in B_0[\bar{q}]$ and yet
$\bar{Q}\in B$. This implies that
$B_0[\bar{q}]\subsetneq  B$,
which is a contradiction.

\section{Maximal Subalgebras That Have No Maximal Subalgebra}
We begin by extending Proposition 1 of [9].
 There the authors showed that if $x\in Z(q)$ with $q$
a sparse Blaschke product, then the algebra
$H^{\infty}_{\supp u_x}$ is a maximal
subalgebra of $H^{\infty}_{\supp u_x}[\bar{q}]$.
Below we show that this is true for a larger class.

\begin{Theorem}
Let $A$ be any Douglas algebra with maximal subalgebra and $x$ be a minimal
support point of
$G_A$. Then $H^{\infty}_{\supp u_x}$ is a maximal subalgebra of
$A_{\supp u_x}$ and $A_{\supp u_x}=H^{\infty}_{\supp u_x}[\bar{q}]$
for some $q\in \Omega(A)$.
\end{Theorem}

{\it Proof:} Let $B_0=H^{\infty}_{\supp u_x}$ and
$B=A_{\supp u_x}$. Suppose
$x$ is a minimal support set for $G_A$. Then we have for
any interpolating Blaschke product $\psi \in \Omega (A)$
with $|\psi (x)|<1$ and any $y\in M(H^{\infty}+C)$
with $\supp u_y
\subsetneq \supp u_x,\,|\psi(y)|=1$. Thus if
$\psi _0\in \Omega(B)$ and $|\psi_0 (x)|<1$
there is a $\psi \in \Omega(A)$ such that $\left.\psi\right|_{\supp u_x}
=\left.\psi_0\right|_{\supp u_x}$.
 This implies that $|\psi _0(y)|=1$ for every such $y$. Hence $x$ is a minimal
support point for $G_B$. Note that this implies that
$M(B)=M(B_0)\setminus E_x$ where $E_x$ is the level
set for $x$. Hence $M(B_0)=M(B)\cup E_x$, so by
Theorem 1 of [10], $B_0$ is a maximal subalgebra of $B$.
Let $q$ be any element in $\Omega (A)$ with $q(x)=0$.
Then $q\in \Omega(B,B_0)$ and we have that $B=B_0[\bar{q}]$.

\begin{Theorem}
A Douglas algebra $A$ has no maximal subalgebra if and only if
$H^{\infty}_{\supp u_x}$ is not a maximal subalgebra of
$A_{\supp u_x}$ for every $x\in G_A$.
\end{Theorem}

{\it Proof:} Suppose $A$ has no maximal subalgebra and let $x\in G_A$.
Since $x$ is not minimal support set of $G_A$ there is a $y\in G_A$ with
$\supp u_y\subsetneq \supp u_x$,
and a $\psi \in \Omega (A)$ such that $|\psi(y)|<1$. Since
$\bar{\psi} \not\in H^{\infty}_{\supp u_y}$, we
can assume that $\psi (y)=0$. Hence
$y\not\in M(A_{\supp u_x})$. By Lemma 4 of [8]
there is a $\psi _0 \in \Omega(A_{\supp u_x})$ such that
$|\psi _0(y)|=1$ and $\psi _0(x)=0$. Then we have
$H^{\infty}_{\supp u_x}\subsetneq 
H^{\infty}_{\supp u_x}[\bar\psi _0]
\subsetneq A_{\supp u_x}$.
So $H^{\infty}_{\supp u_x}$ is not a maximal
subalgebra of $A_{\supp u_x}$.

Suppose that for all $x\in G_A$, $H^{\infty}_{\supp u_x}$
is not a maximal subalgebra of $A_{\supp u_x}$. Then there is
an algebra $B$ with $H^{\infty}_{\supp u_x}
\subsetneq B\subseteq A_{\supp u_x}$.
Thus we can find a $y\in M(H^{\infty}+C)$ such that
$\supp u_y\subsetneq \supp u_x$ and
$y\in M(B)\setminus M(A_{\supp u_x})$. This implies that there is an
interpolating Blaschke product $q$ with $\bar{q}\in B\subset
A_{\supp u_x}$ such
that $|q(y)|=1$ and $|q(x)|<1$. Hence there is a $q_0\in \Omega (A)$ with
$\left.q_0\right|_{\supp u_x}=\left.q\right|_{\supp u_x}$.
So $|q_0(y)|=1$ and $|q_0(x)|<1$. This implies that $x$ is not a minimal
support point of $G_A$ for every $x\in G_A$. So by Theorem 1 of
[10] $A$ has no maximal subalgebra.

\begin{Proposition}
Let $x\in M(H^{\infty})\setminus M(L^{\infty})$.
Then $H^{\infty}_{\supp u_x}$ has no maximal subalgebra.
\end{Proposition}

{\it Proof:} Now $\mathrm{ess}(H^{\infty}_{\supp u_x})=\supp u_x$.
Hence if $y\in G_{H^{\infty}_{\supp u_x}}$, then
$\supp u_y\cap\supp u_x=\emptyset $.
Hence if $A$ is a subalgebra of $H^{\infty}_{\supp u_x}$,
then there is a $y\in M(A)\setminus M(H^{\infty}_{\supp u_x})$
with $\supp u_y\cap \supp u_x=\emptyset$.
Hence $\mathrm{ess}(A)\supseteq \supp u_y\cup
\supp u_x\supsetneq\supp u_x=
\mathrm{ess}(H^{\infty}_{\supp u_x})$.
By Corollary 1, $A$ is not a maximal subalgebra of
$H^{\infty}_{\supp u_x}$.

\begin{Proposition}
Let $A$ be a Douglas algebra that has only one maximal subalgebra $A_0$.
Then $A_0$ has no maximal subalgebra.
\end{Proposition}

{\it Proof:} Suppose there is a subalgebra $B_0\subseteq A_0$ such that
$B_0$ is a maximal subalgebra of $A_0$. Then by Theorem 1 of [10]
there is an $x_0\in G_{A_0}$ such that
\begin{equation}
M(B_0)=M(A_0)\cup E_{x_0}     \label{prop2:1}
\end{equation}
Since $A_0$ is a maximal subalgebra of $A$ there is an $x\in G_A\cap M(B_0)$ such that
\begin{equation}
M(A_0)=M(A)\cup E_x    \label{prop2:2}
\end{equation}

By (1) and (2) we have that
$M(B_0)=M(A)\cup E_x\cup E_{x_0}$. Since $x\in M(A_0)$ we have that
$\supp u_x\neq
\supp u_{x_0}$. Also since $x_0\not\in E_x$ and
$x_0\in G_{A_0}$ implies that $\supp u_{x_0}\not\subset\supp u_x$
(otherwise $x_0\in M(A)$ by (2). We show that $x_0$ is a
minimal support point of
$G_A$, and hence get a contradiction. Let $y\in M(H^{\infty}+C)$
such that $\supp u_y\subsetneq \supp u_{x_0}$.
Since $x_0$ is a minimal support point of  $G_{A_0}$ we have that
$y\in M(A_0)=M(A)\cup E_x$. If $y\in M(A)$ then we are done. So we can
assume $y\in E_x$.
If $y\in E_x$, then $\supp u_y=\supp u_x$, so we have that
$\supp u_x\subsetneq \supp u_{x_0}$.
Since $x\not\in M(A)$ there is an interpolating Blaschke product
$q$ with $\bar{q}\in A$ and such that $q(x)=0$. By Theorem 2
of [11] there is an uncountable set $U$ of $Z(q)$
such that (a) $\supp u_m\subsetneq 
\supp u_{x_0}$ for all $m\in U$ and (b)
$\supp u_m\cap\supp u_k$ for all
$m,k\in U, m\neq k$. By (1) each such $m\in U$ is in $M(A_0)$.
Since
for all $m\in U$ (except if $m=x$) we have $\supp u_x\cap
\supp u_m=\emptyset$,
hence by (2), $m\in M(A)$. But $\bar{q}\in A$ and $U\subset Z(q)\cap M(A).$ This is a contradiction, and we get $y\not\in E_x .$ So $y\in M(A)$ and since $\supp u_{x_0} \neq \supp  u_x$
we have that $x_0$ is a
minimal support point of $G_A$. This is a contradiction. So $A_0$ has no
maximal subalgebra.

 Note that Proposition 1 follows from Proposition 2 if $x$ is a minimal
support point for some
interpolating Blaschke product.

Let $q$ be an interpolating Blaschke product. We consider the algebra
$H^{\infty}_{N(\bar{q})}$. Certainly $H^{\infty}_{N(\bar{q})}$ is
not known to be a maximal subalgebra of any Douglas algebra,
but does have some of the same properties of $H^{\infty}_{\supp u_x}$.
For example

\begin{Proposition}
 The algebra $H^{\infty}_{N(\bar{b})}$ has no maximal subalgebra.
\end{Proposition}

{\it Proof:} 
Set $B=H^{\infty}_{N(\bar{q})}$. Let $x\in G_B$
and suppose $x$ is a minimal support point for $G_B$.
Then if $y\in M(H^{\infty}+C)$ such that $\supp u_y
 \subsetneq   \supp u_x$ we have that $y\in M(B)$.
By [2, page 39] we must have that
$\supp u_y\subseteq N(\bar{q})=\mathrm{ess}(B)$.
Thus we have that $\mathrm{ess}(B)\cap \supp u_x
=N(\bar{q})\cap \supp u_x \neq \emptyset$.
By [14, Theorem 1] $N(\bar{q})=\bigcup_{x\in Z(q)}Q_x$.
So there is an $x_0\in Z(q)$ such that
$\supp u_x \cap Q_{x_0}\neq \emptyset$.
By the definition of $Q_{x_0}$, we have that
$\supp u_x \subset Q_{x_0}$. By [2, page 39] this
implies that $x\in M(B)$, which is a contradiction. So if $x\in G_B$,
then $\supp u_x \cap N(\bar{q})=\emptyset$, which
implies that $x$ is not a minimal support point for $G_B$.
$B$ has no maximal subalgebra.

\section{Minimal Superalgebras of $H^{\infty}_{\supp u_x}$}
We will compute the Bourgain algebras and the Minimal Envelopes
of the Douglas algebra $H^{\infty}_{\supp u_x}$ for any
$x\in M(H^{\infty}+C)$. We have the following

\begin{Theorem}
Let $x\in M(H^{\infty}+C)\setminus M(L^{\infty}$ such that
$|q(x)|<1$ for some interpolating Blaschke product $q$ and set
$B=H^{\infty}_{S_{\supp u_x}}$. Then
\begin{enumerate}
\item[i)] Either $B_b=B$ or $B_b=B[\bar{\psi}]$ for
some interpolating Blaschke product $\psi$.
\item[ii)] Either $B_m=B_b=B$ or $B_m=B[\bar{\psi}]$
for some interpolating Blaschke product $\psi$.
\end{enumerate}
\end{Theorem}
{\it Proof:} We will use Theorem 2 of [4] which
says that for any interpolating Blaschke product $\psi$ with
$\bar{\psi}\in B_b$, we have that the set $Z(\psi)\cap M(B)$
is a finite set and the fact that
\begin{equation}
M(B)=M\left(L^{\infty}\right)\cup \left\{m\in M(H^{\infty}) :
\supp u_m \subseteq \supp u_x \right\} \label{eq:3}
\end{equation}

We claim that if $\psi$ is an interpolating Blaschke product
such that $\bar{\psi} \in B_b$, then $Z(\psi)\cap M(B)\subset E_x$,
the level set of $x$.
Suppose not. Then there is an $x_0\in Z(\psi )\cap M(B)$ such that
$\supp u_{x_0}\subsetneq \supp u_x$.
By Theorem 2 of [11] there is an uncountable set
$\Gamma$ of $Z(\psi)$ such that (a)
$\supp u_{\gamma}\subsetneq \supp u_x$
for all $\gamma \in \Gamma$ and (b)
$\supp u_{x_m}\cap \supp u_{\gamma}=\emptyset$
for all $m,\gamma \in \Gamma,\,m\neq \gamma $.
By (a) and (3) each $\gamma \in M(B)$ and so
$\Gamma \subset Z(\psi)\cap M(B)$. This implies that the set
$Z(\psi)\cap M(B)$ is infinite. This is a contradiction. Hence if
$x_0\in Z(\psi)\cap M(B)$ then
$\supp u_{x_0}= \supp u_{x}$, so we get
$Z(\psi)\cap M(B)\subset E_x$. There are two possibilities
(1) the set $Z(\psi)\cap M(B)=\emptyset$ for which $\bar{\psi}
\in B$, so $B_b\subseteq B$. This gives us the case when $B_b=B$.
(2) If $Z(\psi)\cap M(B)\neq \emptyset$ but finite. Then the
algebra $B[\bar{\psi}]\subseteq B_b$ To show that $B_b=
B[\bar{\psi}]$, let $\psi _0$ be other interpolating Blaschke
product with $\bar \psi _0 \in B_b$. Since both sets
$Z(\psi _0)\cap M(B)$ and $Z(\psi)\cap M(B)$ are contained in
$E_x$, we have that $M(B)\setminus M(B[\bar \psi _0])=
M(B)\setminus M(B[\bar{\psi}])=E_x$. Thus
$M(B[\bar{\psi _0}])=M(B[\bar{\psi}])$ and by
the Chang--Marshall Theorem [1, 16] we have
$B[\bar{\psi _0}]=B[\bar{\psi}]$,
since this is true for all $\psi , \psi _0$ we have by
Theorem C of [12], $B_b=B[\bar{\psi}]$
for any such $\psi$ $B[\bar{\psi}] $ is a minimal
superalgebra of $B$.) This proves  (i).

To prove (ii) let $\bar{\psi}\in B_m$. Then, by Theorem 3 of
[12] there is a finite set
$\{x_1,\dots ,x_n \}\subset Z(\psi )\cap M(B)$ such that
$\{u\in M(B) : |\psi (u)|<1 \} =E_{x_1}\cup \dots \cup E_{x_n}$.
Again we claim that $E_{x_1}=E_{x_2}=\dots =E_{x_n}=E_x$.
Suppose that $E_{x_1}\neq E_{x_2}$. Then
$\supp u_{x_1}\neq \supp u_{x_2}$.
By (3) either
$\supp u_{x_1}\subsetneq \supp u_x$
or
$\supp u_{x_1}\subsetneq \supp u_x$
or both. Let us suppose that
$\supp u_{x_1}\subsetneq \supp u_x$.
Then by Theorem 2 of [11] there is an uncountable
set $\Gamma$ such that $E_{\alpha}\neq E_{\beta}$ for all
$\alpha , \beta \in \Gamma$ and
$\bigcup _{\alpha \in \Gamma} E_{\alpha} \subset \{u \in M(B) :
|\psi (u)|<1 \}$. This contradicts Theorem 3 of [12].
Thus $E_{x_1}=E_{x_2}=\dots =E_{x_n}=E_x$. As before we have that for
$\bar{\psi} \in B_m$, $Z(\psi )\cap M(B) \subset E_x$ and
$B_m=B[\bar{\psi}]$ if $Z(\psi )\cap M(B)\neq \emptyset $.
This proves ii).

\begin{Corollary}
\begin{enumerate}
\item[i)] Let $x\in M(H^{\infty}+C)\setminus M(L^{\infty})$ and
$B=H^{\infty}_{S_{u_x}}$. Then $B\subset B_m$ if and only if $x$
is a minimal support point of $H^{\infty}[\bar{\psi}]$
for some interpolating Blaschke product $\psi $.
\item[ii)] $B=B_b=B_m$ if and only if $x$ is not a minimal support
point of $H^{\infty}[\bar{\psi}]$ for any interpolating
Blaschke product.
\end{enumerate}
\end{Corollary}

Theorem 4 i) has also appeared in [17]. 
\end{document}